# On Mathematical Functions for Theoretical and Experimental Distributions for Shrutis (Micro Tonal Intervals) and Their Application in Hindustani Music

Vishal Midya[1]

[1]*Indian Statistical Institute, Delhi, India*

In this work, exact mathematical functions have been formulated for three important theoretical Shruti (micro tonal interval) distributions, i.e. for Western Compilation, Deval, and Nagoji Row in Hindustani music. A generalized mathematical function for Shrutis has also been formulated. This generalized function shows a very high order of conformity with the experimentally derived Shruti distribution, than those of the theoretical Shruti distributions.

## 1) Introduction

The Indian Musical system is based on modal format, and here the relation between successive notes as well as, the relation between a fixed tonic and any note is of grave importance. Even in the absence of tonic, an isolated note can convey an expression, of course with respect to a memorized tonic or a previous note. Here each note leaves an impression in mind and the idea is revealed through the cumulative effect of successive notes. It should be kept in mind that, the tonic is not fixed, in terms of frequencies, rather can be specified at any position according to the performer. So, in the scale of Indian musical system, relation between successive notes has to be considered.

The notion of Shruti or microtonal intervals has to be given importance. Musicological research on the existence of Shruti and their relevance to the performance or perception has been started long ago. In the ancient system, the swara (musical note) was not merely a sound of fixed pitch position, but the entire tonal range between itself and its previous swara. Although this interval can be divided into infinitesimal parts, it was believed that only a limited number, not exceeding four sounds could be distinctly cognized by the ear in a swara interval. These cognizable sounds are known as Shrutis and the interval, which separated one swara from the next, was measured in terms of Shrutis. The total number of Shrutis was fixed unambiguously at 22 in ancient treatises. Shruti was thought of both as the least audible interval between two sounds, as well as the sounds themselves, which were separated by such an interval. However the two different distributions of Shruti for swara intervals for different gramas (scale types the notes belong to) were reported in the ancient treaties.

Though the awareness of discrepancy between the ancient treatises and contemporary practices are growing at least from the mid 17$^{th}$ century, it was only the latter half of the 20$^{th}$ century that theories of intonation based on empirical research emerged. The western bias for discrete pitches is not only fed by prominence of keyboard instrument; it also has a long history of obsession with ratios. From Pythagoras to Helmholtz, western thinking about melody is marked by a reduction of the continuous tonal space into a limited series of points. Although Indian music uses a similar abstraction by defining 7 notes, 12 semi-tones and 22 Shrutis, it must be stressed that from earlier these may have been considered regions rather than points.

In modern time period, views range from complete negation of practical relevance or even existence of different numbers of them. However, the majority of the authors on Shruti believe that in contemporary music Shrutis exist and form the basis of 12 notes but much of the historical commentaries leave us as baffled as ever. Confusion becomes even more confounded by the unending arithmetical calculations of ratios with little experimental support.

An experimental research on Shrutis has been carried out by A.K.Dutta, R.Sengupta, N.Dey and D.Nag of Scientific Research Department, ITC Sangeet Research Academy and also has been compiled in the book, "Experimental Analysis of Shrutis from the Performances in Hindustani Music". This research shows that, experimentally derived Shruti distribution [1] from the performances of singers in North Indian classical music is different from those of the theoretical Shruti distributions which exist. Though among those, Nagoji Row, Deval and Western Compilation [2] show more or less a good conformity with the experimentally derived Shruti distribution and also indicate the existence of a personal scale, in the form of a personal set of discrete frequencies, embedded in the Indian musical scale based on general Shruti-intervals.

The distribution of Shrutis (micro tonal intervals) in Hindustani musical scale have much complex pattern than that of Western musical scale which follows a logarithmic pattern(Tempered Chromatic scale) in general. Here in Hindustani music the tonic (first note of the musical scale) is not fixed in terms of frequencies rather can be specified at any position according to the performer. Due to these complexities, till now any mathematical function for the theoretical distributions and experimental distribution of Shrutis hasn't been discovered (Western Compilation, Deval, and Nagoji Row are the three most important theoretical Shruti distributions in Hindustani music). So, it has become an age old unsolved problem in computational musicology. Consequently, the lack of any mathematical function for Shruti distributions has hindered the possibility of making any electronic keyboard instrument for Shrutis. Till now most of the Indian Classical musicians use string instruments in performances and don't use any electronic instrument (synthesizer) as all of these are based on Tempered Chromatic scale.

In this context, there comes the need of a generalized Shruti function, which will show a very high order of conformity with the experimentally derived Shruti distribution, than those of the theoretical Shruti distributions, and also meet the purpose of an electronic synthesizer[3], where a performer of Hindustani music on giving his preferable tonic, will get the set of frequencies, which will be used by him in the performance. This work is based on the topics above mentioned.

## 2) Discussion on the Shruti functions

The work starts with the Shruti distribution of western compilation, Deval and Nagoji Row. One function for each of the three theoretical Shruti distributions, as well as their difference in cents from the experimentally derived values of Shrutis is given. Brief discussions on whether any of these three distributions can be used as the generalized mathematical function for Shrutis, are given also below. Lastly there is a discussion on the generalized Shruti equation.

## 2.1) Western Compilation

For Western compilation the function $f(z)$ is defined. This function is got intuitively (At first it was noted that each Shruti ratio of the Western Compilation can be written in terms of 2,3,5,7 i.e. multiplying them together with respect to some integral powers including zero. Then it was also noted that when the 2$^{nd}$ Shruti ratio is divided by the 1$^{st}$ one, the 3$^{rd}$ Shruti ratio is divided by the 2$^{nd}$ one and if we go on dividing like this till the last Shruti ratio, the ratios which we'll get by these continuous divisions will follow a particular pattern. Using these information and with intuition the function $f(z)$ is formed). In this function, the variable $z$ takes integer values, starting from 1,2,3,4,...,to 23 and corresponding to each value of $z$ the function $f(z)$ generates exactly same Shruti value as the Western compilation.

---

[1] Experimentally derived Shruti distribution was got by an objective analysis of **_aalap_** (introductory, slow movement, free of rhythm; unmetered raga introduction and expansion; prelude) from 142 songs performed by 42 eminent musicians and scholars covering 21 different ragas(a melodic concept within certain strict rules. Ragas are specific groups and sequences of notes that organise melody)

[2] Nagoji Row, Deval and Western Compilation are three important Shruti distributions in Hindustani classical music. These were formed empirically.

[3] Till date there is no electronic tool (synthesizer) for Shrutis, due to the lack of function which can generate Shrutis Ratios.

$$f(z) = \left\{ 2^{\left(8\left[\left(\frac{z+2-h(z)}{4}\right)+g_1^2(z)-g_1(z)\right]-4\left[\left(\frac{z-1-h(z)}{2}\right)+g_1(z)\right]-3\left[\left(\frac{z-h(z)}{4}\right)-g_1^2(z)\right]\right)} \right.$$
$$\times 3^{\left(4\left[\left(\frac{z-1-h(z)}{2}\right)+g_1(z)\right]-5\left[\left(\frac{z+2-h(z)}{4}\right)+g_1^2(z)-g_1(z)\right]-\left[\left(\frac{z-h(z)}{4}\right)-g_1^2(z)\right]\right)}$$
$$\left. \times 5^{\left(2\left[\left(\frac{z-h(z)}{4}\right)-g_1^2(z)\right]-\left[\left(\frac{z-1-h(z)}{2}\right)+g_1(z)\right]\right)} \right\}$$

Where:
1) $z = 1(1)23$ And $z$ is a positive integer
2) $g_1(z) = \begin{cases} 0, & z = 1(1)12 \\ -1/2, & z = 13(1)23 \end{cases}$
3) $h(z) = \begin{cases} 1, & when\ z = 13 \\ 0, & othewise \end{cases}$
4) [ ] Denotes greatest integer function

The Shruti values of $f(z)$ and the values of the Shrutis from Western Compilation as well as the difference between the Shruti values of $f(z)$ and Shruti values of experimentally derived Shruti distribution are given below.

**Table for Western Compilation (Table 1)**

| Serial no. | $f(z)$ | Shruti (Western Compilation) | Shruti(Experimental Values) | Difference between W.C. and Shruti exp.(in cents)[4] |
|---|---|---|---|---|
| 01. | 001/001 | 001/001 | 1.000 | 00.00 |
| 02. | 256/243 | 256/243 | 1.048 | 09.06 |
| 03. | 016/015 | 016/015 | 1.065 | 02.71 |
| 04. | 010/009 | 010/009 | 1.117 | 09.15 |
| 05. | 009/008 | 009/008 | 1.134 | 13.79 |
| 06. | 032/027 | 032/027 | 1.189 | 05.56 |
| 07. | 006/005 | 006/005 | 1.208 | 11.50 |
| 08. | 005/004 | 005/004 | 1.227 | 32.15 |
| 09. | 081/064 | 081/064 | 1.267 | 01.88 |
| 10. | 004/003 | 004/003 | 1.328 | 06.94 |
| 11. | 027/020 | 027/020 | 1.349 | 01.28 |
| 12. | 045/032 | 045/032 | 1.399 | 08.95 |
| 13. | 064/045 | 064/045 | 1.429 | 08.23 |
| 14. | 003/002 | 003/002 | 1.506 | 06.91 |
| 15. | 128/081 | 128/081 | 1.579 | 01.37 |
| 16. | 008/005 | 008/005 | 1.604 | 04.32 |
| 17. | 005/003 | 005/003 | 1.630 | 38.51 |
| 18. | 027/016 | 027/016 | 1.682 | 05.65 |

---

[4] **W.C.** – Western Compilation, **exp.**- experimental, **cents**- 1200[th] part of an octave in logarithmic scale

| | | | | |
|---|---|---|---|---|
| 19. | 016/009 | 016/009 | 1.739 | 38.18 |
| 20. | 009/005 | 009/005 | 1.791 | 08.68 |
| 21. | 015/008 | 015/008 | 1.824 | 47.74 |
| 22. | 243/128 | 243/128 | 1.908 | 08.70 |
| 23. | 002/001 | 002/001 | 2.000 | 00.00 |
| Average difference | | | | 12.92 |

The values of the function $f(z)$ are given in the column 2 of table 1. In the 3rd column, the theoretical Shruti distribution of Western Compilation is shown. The function $f(z)$ can be treated as the function for Western Compilation. Now we want to find whether this function can be used as a generalized Shruti function (This generalized Shruti function will be such a function that the average difference between the Shruti values generated by the function and the experimentally obtained Shruti values will be as low as possible). So the difference between the experimentally obtained Shruti values and the Shruti values generated by $f(z)$, are calculated. In the column 4 of table 1, the experimentally obtained Shruti values are given and in the next column, i.e. in the 5th column, the corresponding differences are given (in cents). The average difference comes out to be 12.92 cents. The differences for the 5th, 7th, 8th, 17th, 19th, 21st Shrutis are high (i.e. more than 10 cents). So, the function can't be used as generalized Shruti function, as the average difference must be lower (less than 10 cents) for that.

## 2.2) Deval

Now for Deval Shruti distribution, the function $f_1(z)$ is defined. This function is got by changing some of the dependent variables of $f(z)$. Here also $z$ takes the integer values from 1, 2, 3, 4, 5… to 23 and for each value of $z$ the function $f_1(z)$ gives the Shruti value, which is exactly equal to the corresponding Shruti value of Deval.

$$f_1(z) = \left\{ 2^{\left(8\left[\left(\frac{z+2-h_1(z)}{4}\right)+g_1^2(z)-g_1(z)\right]-4\left[\left(\frac{z-1-h_1(z)}{2}\right)+g_1(z)\right]-3\left[\left(\frac{z-h_1(z)}{4}\right)-g_1^2(z)\right]-10k(z)+4t(z)\right)} \right.$$
$$\times 3^{\left(4\left[\left(\frac{z-1-h_1(z)}{2}\right)+g_1(z)\right]-5\left[\left(\frac{z+2-h_1(z)}{4}\right)+g_1^2(z)-g_1(z)\right]-\left[\left(\frac{z-h_1(z)}{4}\right)-g_1^2(z)\right]+6k(z)-4t(z)\right)}$$
$$\times 5^{\left(2\left[\left(\frac{z-h_1(z)}{4}\right)-g_1^2(z)\right]-\left[\left(\frac{z-1-h_1(z)}{2}\right)+g_1(z)\right]-k(z)+t(z)\right)} \times 7^{(k(z))} \right\}$$

Where:
1) $z = 1(1)23$ And $z$ is a positive integer
2) $g_1(z) = \begin{cases} 0, & z = 1(1)12 \\ -1/2, & z = 13(1)23 \end{cases}$
3) $k(z) = \begin{cases} 1, & z = 2^{\left(1-\left[\frac{n-1}{2}\right]\right)} \times 3^{\left[\frac{n-1}{2}\right]} \times 5^{\left[\frac{n}{2}\right]}, \text{ where } z = 1(1)23 \text{ and } n \in N \\ 0, & \text{otherwise} \end{cases}$
4) $h_1(z) = \begin{cases} 1, & z = 11 \text{ and } z = 13 \\ 0, & \text{otherwise} \end{cases}$
5) $t(z) = \begin{cases} 1, & z = 10 \\ 0, & \text{otherwise} \end{cases}$
6) [ ] Denotes greatest integer function

The Shruti values of $f_1(z)$ and the values of the Shrutis from Deval as well as the difference between the Shruti values of $f_1(z)$ and Shruti values of experimentally derived Shruti distribution are given below.

**Table for Deval (Table 2)**

| Serial no. | $f_1(z)$ | Shruti (Deval) | Shruti(Experimental Values) | Difference between Deval and Shruti exp. (in cents) |
|---|---|---|---|---|
| 01. | 001/001 | 001/001 | 1.000 | 00.00 |
| 02. | 021/020 | 021/020 | 1.048 | 03.30 |
| 03. | 016/015 | 016/015 | 1.065 | 02.71 |
| 04. | 010/009 | 010/009 | 1.117 | 09.15 |
| 05. | 009/008 | 009/008 | 1.134 | 13.79 |
| 06. | 032/027 | 032/027 | 1.189 | 05.56 |
| 07. | 006/005 | 006/005 | 1.208 | 11.50 |
| 08. | 005/004 | 005/004 | 1.227 | 32.15 |
| 09. | 081/064 | 081/064 | 1.267 | 01.88 |
| 10. | 021/016 | 021/016 | 1.328 | 20.33 |
| 11. | 004/003 | 004/003 | 1.349 | 20.22 |
| 12. | 045/032 | 045/032 | 1.399 | 08.95 |
| 13. | 064/045 | 064/045 | 1.429 | 08.23 |
| 14. | 003/002 | 003/002 | 1.506 | 06.91 |
| 15. | 063/040 | 063/040 | 1.579 | 04.39 |
| 16. | 008/005 | 008/005 | 1.604 | 04.32 |
| 17. | 005/003 | 005/003 | 1.630 | 38.51 |
| 18. | 027/016 | 027/016 | 1.682 | 05.65 |
| 19. | 016/009 | 016/009 | 1.739 | 38.18 |
| 20. | 009/005 | 009/005 | 1.791 | 08.68 |
| 21. | 015/008 | 015/008 | 1.824 | 47.74 |
| 22. | 243/128 | 243/128 | 1.908 | 08.70 |
| 23. | 002/001 | 002/001 | 2.000 | 00.00 |
| **Average difference** | | | | **14.33** |

In the column 2 and 3 of table 2, the values of $f_1(z)$ and the corresponding Shruti values of Deval are given respectively. As $f_1(z)$ gives the same Shruti values as Deval, we can use the function $f_1(z)$ as the function for Deval Shruti distribution. Now to find if the function $f_1(z)$ can be used as the generalized Shruti function, the differences between the experimentally obtained Shruti values and the corresponding Shruti values generated by $f_1(z)$ are calculated (in cents). In the column 4 of table 2, the experimentally obtained Shruti values are given and in the next column the corresponding differences are shown. The average difference comes out to be 14.33 cents, which is obviously greater than the average difference of Western Compilation. Here, the differences for the 5th, 7th, 8th, 10th, 11th, 17th, 19th, 21st Shrutis are high (more than 10 cents) and hence the average difference is high. So, $f_1(z)$ can't be used as the generalized Shruti function.

## 2.3) Nagoji Row

For Nagoji Row, the function $f_2(z)$ is defined. This function is also got by changing some dependent variables of $f(z)$. As previous, here also z takes integer values, starting from 1,2,3,4, ...,to 23 and corresponding to each integer value of z, the function $f_2(z)$ generates Shruti value, which is exactly equal to the corresponding Shruti value of Nagoji Row.

$$f_2(z) = \left\{ 2^{\left(8\left[\left(\frac{z+2}{4}\right)+g^2(z)-g(z)\right]-4\left[\left(\frac{z-1}{2}\right)+g(z)\right]-3\left[\left(\frac{z}{4}\right)-g^2(z)\right]-11p(z)\right)} \right.$$

$$\times\, 3^{\left(4\left[\left(\frac{z-1}{2}\right)+g(z)\right]-5\left[\left(\frac{z+2}{4}\right)+g^2(z)-g(z)\right]-\left[\left(\frac{z}{4}\right)-g^2(z)\right]+4p(z)\right)}$$

$$\left. \times\, 5^{\left(2\left[\left(\frac{z}{4}\right)-g^2(z)\right]-\left[\left(\frac{z-1}{2}\right)+g(z)\right]+2p(z)\right)} \right\}$$

Where:

1) $z = 1(1)23$ And $z$ is a positive integer

2) $g(z) = \begin{cases} 0, & z = 1(1)13 \\ -1/2, & z = 14(1)23 \end{cases}$

3) $p(z) = \begin{cases} 1, & \text{for } \left[\frac{z}{7}\right] = 0 \text{ or even positive integer}; \\ & \text{when } z = 2 + 7\left[\frac{n+1}{3}\right] + 4\left[\frac{n}{3}\right] + 2\left[\frac{n-1}{3}\right], \text{where } n \in N \text{ and } z = 1(1)23 \\ -1, & \text{for } \left[\frac{z}{7}\right] = \text{odd positive integer}; \\ & \text{when } z = 2 + 7\left[\frac{n+1}{3}\right] + 4\left[\frac{n}{3}\right] + 2\left[\frac{n-1}{3}\right], \text{where } n \in N \text{ and } z = 1(1)23 \end{cases}$

4) [ ] Denotes greatest integer function

The Shruti values of $f_2(z)$ and the values of the Shrutis from Nagoji Row as well as the difference between the Shruti values of $f_2(z)$ and Shruti values of experimentally derived Shruti distribution are given below.

**Table for Nagoji Row (Table 3)**

| Serial no. | $f_2(z)$ | Shruti (Nagoji Row) | Shruti(Experimental Values) | Difference between N.R. and Shruti exp.(in cents)[5] |
|---|---|---|---|---|
| 01. | 001/001 | 001/001 | 1.000 | 00.00 |
| 02. | 025/024 | 025/024 | 1.048 | 10.49 |
| 03. | 016/015 | 016/015 | 1.065 | 02.71 |
| 04. | 010/009 | 010/009 | 1.117 | 09.15 |
| 05. | 009/008 | 009/008 | 1.134 | 13.79 |
| 06. | 032/027 | 032/027 | 1.189 | 05.56 |
| 07. | 006/005 | 006/005 | 1.208 | 11.50 |
| 08. | 005/004 | 005/004 | 1.227 | 32.15 |
| 09. | 032/025 | 032/025 | 1.267 | 17.67 |
| 10. | 004/003 | 004/003 | 1.328 | 06.94 |
| 11. | 027/020 | 027/020 | 1.349 | 01.28 |
| 12. | 045/032 | 045/032 | 1.399 | 08.95 |
| 13. | 036/025 | 036/025 | 1.429 | 13.28 |
| 14. | 003/002 | 003/002 | 1.506 | 06.91 |
| 15. | 025/016 | 025/016 | 1.579 | 18.19 |
| 16. | 008/005 | 008/005 | 1.604 | 04.32 |
| 17. | 005/003 | 005/003 | 1.630 | 38.51 |
| 18. | 027/016 | 027/016 | 1.682 | 05.65 |

---

[5] **N.R.** – Nagoji Row

| | | | | |
|---|---|---|---|---|
| 19. | 016/009 | 016/009 | 1.739 | 38.18 |
| 20. | 009/005 | 009/005 | 1.791 | 08.68 |
| 21. | 015/008 | 015/008 | 1.824 | 47.74 |
| 22. | 048/025 | 048/025 | 1.908 | 10.85 |
| 23. | 002/001 | 002/001 | 2.000 | 00.00 |
| **Average difference** | | | | **14.88** |

In the 3$^{rd}$ column of the table 3, the theoretical Shruti distribution of Nagoji Row is shown. The values of the function $f_2(z)$ are given in the column 2 of table 3. The function $f_2(z)$ can be treated as the function for Nagoji Row as the Shruti values generated by $f_2(z)$ are same as Nagoji Row. Now to find whether this function can be used as a generalized Shruti function, the difference between the experimentally obtained Shruti values and the Shruti values generated by $f_2(z)$ are calculated. In the column 4 of table 3, the experimentally obtained Shruti values are given and in the next column i.e. in the 5$^{th}$ column, the corresponding differences are given (in cents). The average difference comes out to be 14.88 cents, which is even greater than the average difference of Deval. The differences for the 2$^{nd}$ ,5$^{th}$ ,7$^{th}$ ,8$^{th}$,9$^{th}$ ,13$^{th}$ ,15$^{th}$ ,17$^{th}$ ,19$^{th}$ ,21$^{st}$ ,22$^{nd}$ Shrutis are high (i.e. more than 10 cents) and hence the average difference is high. So, $f_2(z)$ also can't be used as the generalized Shruti function.

## 2.4) Generalized Shruti Function

Keeping in mind the fact, that the generalized Shruti function will be such that the average difference between the Shruti values generated by the function and the experimentally obtained Shruti values will be as low as possible, the function $F(z)$ is defined. This function is actually got by taking an average of the functions $f_1(z)$ and $f_2(z)$, and altering some of the dependent variables of the both the functions.

$$F(z) = \{ 2^{\left( 8\left[\left(\frac{z+2-\Psi(z)}{4}\right)+g^2(z)-g(z)\right]-4\left[\left(\frac{z-\Psi(z)-1}{2}\right)+g(z)\right]-3\left[\left(\frac{z-\Psi(z)}{4}\right)-g^2(z)\right]-11\alpha(z)+6\beta(z)-7\gamma(z)-1 \right)}$$
$$\times 3^{\left( 4\left[\left(\frac{z-\Psi(z)-1}{2}\right)+g(z)\right]-5\left[\left(\frac{z+2-\Psi(z)}{4}\right)+g^2(z)-g(z)\right]-\left[\left(\frac{z-\Psi(z)}{4}\right)-g^2(z)\right]+4\alpha(z)-2\beta(z)-3\gamma(z) \right)}$$
$$\times 5^{\left( 2\left[\left(\frac{z-\Psi(z)}{4}\right)-g^2(z)\right]-\left[\left(\frac{z-\Psi(z)-1}{2}\right)+g(z)\right]+2\alpha(z)-\gamma(z) \right)} \times 7^{(-\beta(z))}$$
$$+ 2^{\left( 8\left[\left(\frac{z+2}{4}\right)+g^2(z)-g(z)\right]-4\left[\left(\frac{z-1}{2}\right)+g(z)\right]-3\left[\left(\frac{z}{4}\right)-g^2(z)\right]-1 \right)}$$
$$\times 3^{\left( 4\left[\left(\frac{z-1}{2}\right)+g(z)\right]-5\left[\left(\frac{z+2}{4}\right)+g^2(z)-g(z)\right]-\left[\left(\frac{z}{4}\right)-g^2(z)\right] \right)} \times 5^{\left( 2\left[\left(\frac{z}{4}\right)-g^2(z)\right]-\left[\left(\frac{z-1}{2}\right)+g(z)\right] \right)} \}$$

<u>Where:</u>
1) $z = 1(1)23$ And $z$ is a positive integer
2) $g(z) = \begin{cases} 0, & z = 1(1)13 \\ -1/2, & z = 14(1)23 \end{cases}$
3) $\beta(z) = \begin{cases} 1, & when\ z = 5 \\ 0, & othewise \end{cases}$
4) $\alpha(z) = \begin{cases} 1, & for\ \left[\frac{z}{7}\right] = 0\ or\ even\ positive\ integer, \\ & when\ z = 20n - 18, n\epsilon N\ and\ z = 1(1)23 \\ -1, & for\ \left[\frac{z}{7}\right] = odd\ positive\ integer \\ & when\ z = 20n - 18, n\epsilon N\ and\ z = 1(1)23 \\ 0, & otherwise \end{cases}$

5) $\Psi(z) = \begin{cases} 1, & z = 8 + 9\left[\frac{n+1}{3}\right] + 2\left[\frac{n}{3}\right] + \left[\frac{n-1}{3}\right], where\ z = 1(1)23\ and\ n \in N \\ 0, & otherwise \end{cases}$

6) $\gamma(z) = \begin{cases} 1, & when\ z = 21 \\ 0, & othewise \end{cases}$

7) [ ] Denotes greatest integer function

The values of $F(z)$ along with the experimentally derived values of shrutis, as well as their difference in cents are given below:

**Table for $F(z)$ (Table 4)**

| Serial no. | $F(z)$ | Shruti(Experimental Values) | Difference between $F(z)$ Shruti exp.(in cents) |
|---|---|---|---|
| 01. | 1.000 | 1.000 | 0.00 |
| 02. | 1.048 | 1.048 | 0.00 |
| 03. | 1.067 | 1.065 | 3.25 |
| 04. | 1.111 | 1.117 | 9.32 |
| 05. | 1.134 | 1.134 | 0.00 |
| 06. | 1.185 | 1.189 | 5.83 |
| 07. | 1.200 | 1.208 | 11.5 |
| 08. | 1.225 | 1.227 | 2.82 |
| 09. | 1.266 | 1.267 | 1.37 |
| 10. | 1.333 | 1.328 | 6.51 |
| 11. | 1.350 | 1.349 | 1.28 |
| 12. | 1.406 | 1.399 | 8.64 |
| 13. | 1.424 | 1.429 | 6.07 |
| 14. | 1.500 | 1.506 | 6.91 |
| 15. | 1.580 | 1.579 | 1.10 |
| 16. | 1.600 | 1.604 | 4.32 |
| 17. | 1.633 | 1.630 | 3.18 |
| 18. | 1.688 | 1.682 | 6.16 |
| 19. | 1.733 | 1.739 | 5.98 |
| 20. | 1.789 | 1.791 | 1.93 |
| 21. | 1.826 | 1.824 | 1.90 |
| 22. | 1.909 | 1.908 | 0.91 |
| 23. | 2.000 | 2.000 | 0.00 |
| **Average difference** | | | **4.04** |

In the column 2 of the table 4, the Shruti values generated by $F(z)$ are given. In the next column, the experimentally obtained Shruti values are shown. In the column 4 the differences (in cents) between the Shruti values generated by $F(z)$ and the experimentally obtained Shruti values are also shown. The average difference comes out to be 4.04 cents. In the table 5, average differences of the three Shruti distributions, namely Western Compilation, Deval, and Nagoji Row and $F(z)$ are given. As the average difference for $F(z)$ is very low, the function can be considered as the generalized Shruti equation.

**Average difference of Shruti distributions from experimental Shruti values: (Table 5)**

|  | W.C. $(f(z))$ | Deval $(f_1(z))$ | N.R. $(f_2(z))$ | $F(z)$ |
|---|---|---|---|---|
| **Average differences (in cents)** | 12.92 | 14.33 | 14.88 | 04.04 |

From the above table it is clear that the function $F(z)$ has the lowest difference from the experimentally derived Shruti values, among all other Shruti distributions. And a difference of **04.04 Cents** can be easily considered very small. So, the function $F(z)$ meets the condition of being the generalized mathematical Shruti function, as it is highly close to the experimentally obtained Shruti distribution.

## 3) Experimental verification

To experimentally verify, whether the Shruti ratios generated by $F(z)$ are the preferred Shruti positions for the performers in Hindustani music, Shruti ratios of the aalap portion for raga Bhairav and raga Darbari Kannada, sung by 27 different singers have been collected (some of the singers have sung both the alaap portions of the two ragas)[6]. For each singer at first the differences between his/her sung Shruti ratios and the corresponding Shruti ratios generated by $F(z)$ are calculated and for each of them the average difference are calculated too and then the overall average difference is found out. This computation has been done for both the ragas.

**Table for the average differences between sung Shruti ratios and Shruti ratios of $F(z)$ for different singers for Raga Bhairav: (Table 6)**

| Serial no. | Singers | Average difference (in Cents) |
|---|---|---|
| 01. | $S_{01}$ | 09.343 |
| 02. | $S_{02}$ | 05.107 |
| 03. | $S_{03}$ | 17.564 |
| 04. | $S_{04}$ | 07.383 |
| 05. | $S_{05}$ | 13.055 |
| 06. | $S_{06}$ | 08.735 |
| 07. | $S_{07}$ | 06.735 |
| 08. | $S_{08}$ | 12.521 |
| 09. | $S_{09}$ | 14.603 |
| 10. | $S_{10}$ | 10.973 |
| 11. | $S_{11}$ | 07.252 |
| 12. | $S_{12}$ | 07.222 |
| 13. | $S_{13}$ | 09.664 |
| 14. | $S_{14}$ | 10.667 |
| 15. | $S_{15}$ | 07.392 |
| 16. | $S_{16}$ | 08.522 |
| 17. | $S_{17}$ | 06.908 |
| 18. | $S_{18}$ | 12.045 |
| 19. | $S_{19}$ | 10.594 |
| 20. | $S_{20}$ | 13.051 |
| 21. | $S_{21}$ | 06.653 |
| 22. | $S_{22}$ | 03.157 |
| 23. | $S_{23}$ | 04.648 |
| 24. | $S_{24}$ | 08.438 |
| 25. | $S_{25}$ | 07.496 |
| 26. | $S_{26}$ | 12.337 |

---

[6] Data on Shruti ratios of 'Aalap' portion of Raga Bhairav and Raga Darbari Kannada - from 'Experimental Analysis of Shrutis from the Performances in Hindustani Music', Scientific Research Department, ITC Sangeet Research Academy, Kolkata, India.

| 27. | $S_{27}$ | 09.999 |
|---|---|---|
| **Grand Average** | | **09.336** |

Where, $S_n, n = 1(1)27$ denotes the $n^{th}$ singer.

**Table for the average differences between sung Shruti ratios and Shruti ratios of $F(z)$ for different singers for Raga Darbari Kannada: (Table 7)**

| Serial no. | Singers | Average difference (in Cents) |
|---|---|---|
| 01. | $S_{01}'$ | 09.004 |
| 02. | $S_{02}'$ | 09.616 |
| 03. | $S_{03}'$ | 09.318 |
| 04. | $S_{04}'$ | 04.060 |
| 05. | $S_{05}'$ | 08.561 |
| 06. | $S_{06}'$ | 07.947 |
| 07. | $S_{07}'$ | 02.660 |
| 08. | $S_{08}'$ | 12.388 |
| 09. | $S_{09}'$ | 15.446 |
| 10. | $S_{10}'$ | 05.267 |
| 11. | $S_{11}'$ | 11.242 |
| 12. | $S_{12}'$ | 11.955 |
| 13. | $S_{13}'$ | 16.690 |
| 14. | $S_{14}'$ | 07.842 |
| 15. | $S_{15}'$ | 08.222 |
| 16. | $S_{16}'$ | 07.713 |
| 17. | $S_{17}'$ | 08.533 |
| 18. | $S_{18}'$ | 12.417 |
| 19. | $S_{19}'$ | 16.311 |
| 20. | $S_{20}'$ | 12.006 |
| 21. | $S_{21}'$ | 08.028 |
| 22. | $S_{22}'$ | 07.872 |
| 23. | $S_{23}'$ | 07.399 |
| 24. | $S_{24}'$ | 13.117 |
| 25. | $S_{25}'$ | 08.490 |
| 26. | $S_{26}'$ | 09.767 |
| 27. | $S_{27}'$ | 08.645 |
| **Grand Average** | | **09.649** |

Where, $S_n', n = 1(1)27$ denotes the $n^{th}$ singer.

In the table 6 for raga Bhairav, the average difference of each singer is given and in the last row the overall average difference is shown. It comes out to be 9.336 cents. In table 7, for raga Darbari Kannada also, the average difference for each singer is given and in the last row, the overall average is shown. It is 9.649 cents. It should be noted that both the average differences for the two ragas can be considered low. Though these two overall average differences are larger than that of 4.04 cents when the average difference is calculated between the experimentally obtained Shruti values and the Shruti values generated by $F(z)$. The reason may be that the singers chosen for this experimental verification are from different Gharana (sub tradition or school) and hence, the usage or preferences of Shrutis are different. But as a whole, the function $F(z)$ can be considered as the generalized Shruti equation.

# 4) Conclusion and Further Research

The average difference between the Shruti ratios of $F(z)$ and the experimentally obtained values of Shrutis is $4.04$ cents, and that of between the Shruti ratios of $F(z)$ and the sung Shruti ratios of singers of Raga Bhairav and Raga Darbari Kannada are $9.336\ cents$ and $9.649\ cents$. These differences can be considered reasonably low. So, evidently the function $F(z)$ can be considered as the generalized Shruti equation. These functions are the solutions to the unsolved problem of "formulation of mathematical functions for Theoretical Shruti (micro tonal intervals) distributions in Hindustani music" in computational musicology.

To conclude, it should be stated that, this function $F(z)$ is not only of theoretical significance, but also has deep practical use. This function might be helpful for the development of an electronic synthesizer, based on Shrutis, where a performer of Hindustani classical music will get the set of frequencies which will be used by him in performance on the basis of his chosen tonic and the synthesizer may provide him with all the required frequencies with the help of the function. Further research can be done on developing Generalized Shruti equations for which the average difference and also the no. of variables used will be less than $F(z)$.

Also there is a scope of doing research on developing electronic synthesizer for Shrutis using this function. This research opens up a totally new area of study in computational musicology .These expressions for Shruti distributions will help to form mathematical models of not only Hindustani Ragas but also of any piece of music which is based on Shrutis.

# 5) Acknowledgement

I would like to thank my advisor Dr. Ranjan Sengupta, Scientific Advisor, C.V. Raman Centre for Physics and Music, Jadavpur University, Kolkata, India, for his untiring support and numerous technical discussions.